\documentclass[11pt,reqno]{amsart} %
\usepackage{amsmath,amscd}
\usepackage{amsthm}
\usepackage{amssymb}
\usepackage{amsfonts}
\usepackage{latexsym}
\usepackage{url}
\usepackage{latexsym}
\usepackage{graphicx,color,subfig}
\usepackage{enumerate,fancybox}
\usepackage[colorlinks=true, citecolor=blue]{hyperref} % to make \ref and \eqref appear as a colored link

\newcommand{\sect}[1]{\section{#1}}

\newtheorem{thh}{Theorem}[section]

\newtheorem{lem}[thh]{Lemma}

\newtheorem{example}[thh]{Example}

\newtheorem{theorem}{Theorem}[section]
\newtheorem{proposition}{Proposition}[section]
\newtheorem{lemma}{Lemma}[section]

\newtheorem{remark}{Remark}[section]

\def\sep{\;\vrule\;}
\def\proofs#1. {\par
                      \ifdim\lastskip<15pt
                      \removelastskip\penalty-200
                      \vskip15pt plus3pt minus3pt
                      \fi
                       {\def\a{#1}
                       \ifx\a\empty
                       {\noindent\bf Proof.}
                       \else
                       {\noindent\bf Proof of #1.}
                       \fi}\enspace}
\def\restr#1{\,\vrule\,\lower1.75ex\hbox{$#1$}}

\def\endproof{\hfill\hspace{-6pt}\rule[-14pt]{6pt}{6pt}
\vskip22pt plus3pt minus 3pt}

\def\be{\begin{equation}}
\def\ee{\end{equation}}
\def\bea{\begin{eqnarray}}
\def\eea{\end{eqnarray}}
\def\bean{\begin{eqnarray*}}
\def\eean{\end{eqnarray*}}

\def\a{\alpha}

\def\d{\delta}

\def\e{\varepsilon}
\def\f{\varphi}

\def\g{\gamma}
\def\G{\Gamma}
\def\i{\infty}

\def\l{\lambda}

\def\s{\sigma}

\def\c{{\rm cap}}

\def\setm{\setminus}

\def\ov{\overline}

\def\c{{\rm cap}}

\font\tenopen = cmbx10
\font\sevenopen = cmbx7
\font\fiveopen = cmbx5
\newfam\openfam
\def\open{\fam\openfam\tenopen}
\textfont\openfam = \tenopen
\scriptfont\openfam = \sevenopen
\scriptscriptfont\openfam = \fiveopen

\def\C{{\open C}}

\newcommand{\capGm}{{\mathrm{cap}}(\Gamma)}

\def\Gbar{\overline{G}}

\DeclareMathOperator\dist{dist}

\input colordvi

\newcommand{\conj}[1]{\overline{#1}}

\newcommand{\jeqN}{j=1,2,\ldots,m}

\input colordvi

\author{E. B. Saff, H. Stahl\dag, N. Stylianopoulos and V. Totik}

\title[Area-type measures and recovery]{Orthogonal polynomials for area-type measures and image recovery}
\thanks{
{\it Acknowledgements.}
The first author was partially supported by the
U.S.\ National Science Foundation grant DMS-1109266.
The third author was supported by the University of Cyprus grant 3/311-21027.
The fourth author was supported by the U.S.\ National Science Foundation grant DMS-1265375.
All authors are indebted to the Mathematical Research Institute at
Oberwolfach, Germany, which provided exceptional working conditions during a Research in Pairs
workshop in 2011, when this paper was conceived.}
\date{\today}

\begin{document}
\maketitle

\begin{abstract}
Let $G$ be a finite union of disjoint and bounded Jordan domains in the
complex plane, let $\mathcal{K}$ be a compact subset of $G$ and consider the set
$G^\star$ obtained from $G$ by removing $\mathcal{K}$; i.e., $G^\star:=G\setminus \mathcal{K}$.
We refer to $G$ as an archipelago and $G^\star$ as an archipelago with lakes.
Denote by  $\{p_n(G,z)\}_{n=0}^\infty$ and $\{p_n(G^\star,z)\}_{n=0}^\infty$,
 the sequences of the Bergman polynomials
associated with  $G$ and $G^\star$, respectively; that is, the orthonormal polynomials
with respect to the area measure on $G$ and $G^\star$. The purpose of the
paper is to show that $p_n(G,z)$ and $p_n(G^\star,z)$ have comparable asymptotic properties,
thereby demonstrating that the asymptotic properties of the Bergman polynomials for $G^\star$ are
determined by the boundary of $G$. As a consequence we can analyze
certain asymptotic properties of $p_n(G^\star,z)$  by
using the corresponding results for $p_n(G,z)$, which were obtained in a recent work
by B.\ Gustafsson, M.\ Putinar, and two of the present authors. The results lead
to a reconstruction algorithm for recovering the shape of an archipelago
with lakes from a partial set of its complex moments.
\end{abstract}

\maketitle
%\tableofcontents
%\allowdisplaybreaks

\begin{figure}[h]
\begin{center}
\includegraphics[scale=0.80]{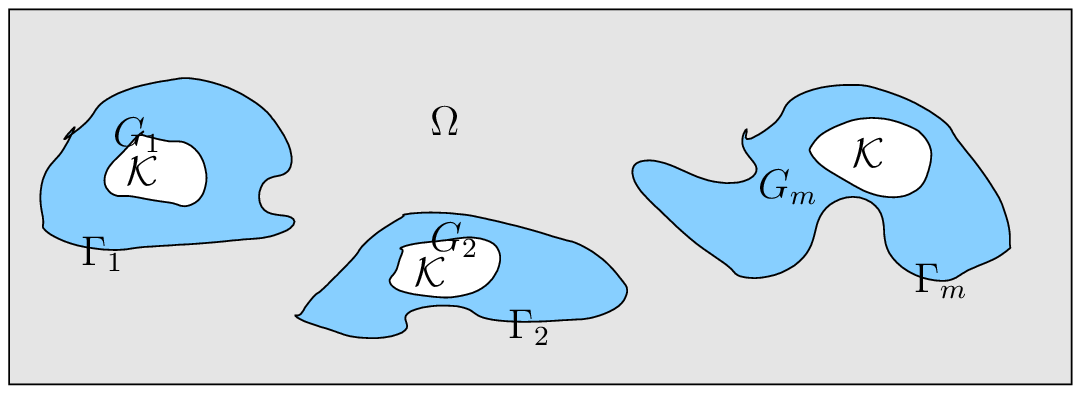}
\caption{}
\label{fig:Archi-lakes}
\end{center}
\end{figure}

\sect{Introduction}\label{sec:intro}
%\subsection{Bergman spaces and polynomials}\label{subsec:bergman}
Let $G:=\cup_{j=1}^m G_j$ be a finite union of bounded Jordan domains $G_j$, $j=1,\ldots,m$, in the
complex plane $\C$, with pairwise disjoint closures, let $\mathcal{K}$ be a compact subset of $G$ and
consider the set
$G^\star$ obtained from $G$ by removing $\mathcal{K}$, i.e., $G^\star:=G\setminus \mathcal{K}$.  Set $\Gamma_j:=\partial G_j$ for the  respective
boundaries and let $\Gamma:=\cup_{j=1}^m \Gamma_j$ denote the boundary of $G$. For later use we introduce also the
(unbounded) complement $\Omega$ of $\overline{G}$ with respect to $\overline{\C}$, i.e.,
$\Omega:=\overline{\C}\setminus\overline{G}$; see Figure~\ref{fig:Archi-lakes}. Note that
$\Gamma=\partial G=\partial\Omega$.
We call $G$ an \emph{archipelago} and  $G^\star$ an \emph{archipelago with lakes}.

Let $\{p_n(G,z)\}_{n=0}^\infty$ denote the sequence of {\em Bergman polynomials} associated with $G$. This is defined as the
unique sequence of polynomials
$$
p_n(G,z) = \gamma_n(G) z^n+ \cdots, \quad \gamma_n(G)>0,\quad
n=0,1,2,\ldots,
$$
that are orthonormal with respect to the inner product
\begin{equation}\label{ee:inprG}
\langle f,g\rangle_G := \int_G f(z) \overline{g(z)} dA(z),
\end{equation}
where $dA$ stands for the differential of the area  measure. We use $L^2(G)$ to denote the associated Lebesgue space with
norm $\|f\|_{L^2(G)}:=\langle f,f\rangle_G^{1/2}$.

The corresponding monic polynomials
$p_n(G,z)/\gamma_n(G)$, can be equivalently defined by the extremal property
\[
\left\|\frac{1}{\gamma_n(G)}{p_n(G,\cdot)}\right\|_{L^2(G)}:=
\min_{z^n+\cdots}\|z^n+\cdots\|_{L^2(G)}.
\]
Thus,
\begin{equation}\label{eq:minimal2}
\frac{1}{\gamma_n(G)}=\min_{z^n+\cdots}\|z^n+\cdots\|_{L^2(G)}.
\end{equation}
A related extremal problem leads to the sequence $\{\lambda_n(G,z)\}_{n=1}^\infty$ of the so-called \emph{Christoffel functions}  associated with the area measure on $G$. These are defined, for any $z\in\C$, by
\begin{equation}\label{eq:Chrfun-def}
\lambda_n(G,z):=\inf\{\|P\|_{L^2(G)}^2,\, P\in\mathbb{P}_n\mbox{ with }  P(z)=1\},
\end{equation}
where $\mathbb{P}_n$ stands for the space of complex polynomials of degree up to $n$.
Using the Cauchy-Schwarz inequality it is easy to verify
(see, e.g., \cite[Section~3]{To05}) that
\begin{equation}\label{eq:Chrfun-pro}
\frac{1}{\lambda_n(G,z)}=\sum_{k=0}^n|p_k(G,z)|^2,\quad z\in\C.
\end{equation}
Clearly, $\lambda_n(G,z)$ is the inverse of the diagonal of the kernel polynomial
\begin{equation}\label{eq:kernel-poly-def}
K^G_n(z,\zeta):=\sum_{k=0}^n \overline{p_k(G,\zeta)} p_k(G,z).
\end{equation}

We use $L_a^2(G)$ to denote the \textit{Bergman space} associated with $G$ and the inner product (\ref{ee:inprG}), i.e.,
\[
{L}_a^2(G):=\left\{ f\ \mathrm{analytic\ in}\ G\ \mathrm{and}\
\|f\|_{{L}^2(G)}<\infty \right\},
\]
and note that $L_a^2(G)$ is a Hilbert space that possesses a reproducing kernel, which we denote by $K_G(z,\zeta)$. That is, $K_G(z,\zeta)$ is the unique function $K_G(z,\zeta):G\times G\to\C$ such that $K_G(\cdot,\zeta)\in {L}_a^2(G)$,
for all $\zeta\in G$, with the reproducing property
\begin{equation}\label{eq:rep_on_L2aG}
f(\zeta)=\langle f,K_G(\cdot,\zeta)\rangle_G, \quad \forall\,  f\in {L}_a^2(G).
\end{equation}
In particular, for any $z\in G$,
\begin{equation}\label{eq:KGzz}
K_G(z,z)=\|K_G(\cdot,z)\|^2_{{L}^2(G)}>0,
\end{equation}
which, in view of the reproducing property and the Cauchy-Schwarz inequality, yields the characterization
\begin{equation}\label{eq:KG(z,z)inv}
\frac{1}{K_G(z,z)}=\inf\{\|f\|_{L^2(G)}^2,\, f\in L_a^2(G)\mbox{ with }  f(z)=1\},
\end{equation}
cf. (\ref{eq:Chrfun-def})--(\ref{eq:kernel-poly-def}).
Furthermore, due to the same property and the completeness of polynomials in $L_a^2(G)$
(see, e.g., \cite[Lemma 3.3]{GPSS}), the kernel $K_G(z,\zeta)$ is given, for any $\zeta\in G$, in terms of the Bergman
polynomials by
\begin{equation}\label{KG}
K_G(z,\zeta)=\sum_{n=0}^\infty \overline{p_n(G,\zeta)}p_n(G,z),
\end{equation}
locally uniformly with respect to $z\in G$.

Consider now the {Bergman spaces} ${L}_a^2(G_j)$, $\jeqN$, associated with the components $G_j$,
\[
{L}_a^2(G_j):=\left\{ f\ \mathrm{analytic\ in}\ G_j\ \mathrm{and}\
\|f\|_{{L}^2(G_j)}<\infty \right\},
\]
and let $K_{G_j}(z,\zeta)$ denote their respective reproducing kernels. Then it is straightforward to
verify using the uniqueness property of $K_G(\cdot,\zeta)$ the following relation
\begin{equation}\label{eq:KK1}
 K_G(z,\zeta)=\left\{
\begin{array}{cl}
K_{G_j}(z,\zeta)  &\mathrm{if}\ z, \zeta\in G_j,\,\,j=1,\ldots,m,\\
 0            &\mathrm{if}\ z\in G_j,\ \zeta\in G_k,\ j\neq k.
\end{array}
\right.
\end{equation}

This relation leads to expressing $K_G(z,\zeta)$ in terms of
conformal mappings $\varphi_j:G_j\to\mathbb{D}$, $\jeqN$.
This is so because, as it is well-known (see
e.g.\ \cite[p.\ 33]{Gabook87}), for $z,\zeta\in G_j$,
\[
K_{G_j}(z,\zeta)=\frac{\varphi'_j(z)\conj{\varphi'_j(\zeta)}}
{\pi\,\left[1-\varphi_j(z)\conj{\varphi_j(\zeta)}\right]^2}.
\]

For  $G^\star:=G\setminus \mathcal{K}$, we likewise define $\langle f,g\rangle_{G^\star}$, the norm $\|f\|_{L^2(G^\star)}$, the Bergman space ${L}_a^2(G^\star)$ along with its reproducing kernel $K_{G^\star}(z,\zeta):G^\star\times G^\star\to\C$ and associated orthonormal polynomials
$$
p_n(G^\star,z) = \gamma_n(G^\star) z^n+ \cdots, \quad \gamma_n(G^\star)>0,\quad
n=0,1,2,\ldots,
$$
as well as the associated  Christoffel functions $\lambda_n^\star(G,z)$ and polynomial kernel functions $ K^{G^\star}_n(z,\zeta).$ It is important to
note, however, that the analogue of \eqref{KG} with $G$ replaced by $G^\star$ does not hold because the polynomials $\{p_n(G^\star,z)\}_{n=0}^{\infty}$
are not complete in ${L}_a^2(G^\star).$

Since $G^\star\subset G$, it is  readily verified that the following two comparison principles hold:
\begin{equation}\label{eq:comp-pri1}
\lambda_n(G^\star,z)\le \lambda_n(G,z),\quad z\in\C,
\end{equation}
and
\begin{equation}\label{eq:comp-pri2}
K_{G}(z,z)\le K_{G^\star}(z,z),\quad z\in G^\star.
\end{equation}

The paper is organized as follows. In the next three sections we prove that holes inside the domains have little influence
on the external asymptotics (a fact anticipated in \cite[Section 3]{Saff1}). Then, in the last section we use this to modify the recent domain recovery algorithm
from \cite{GPSS} to the case when one has no \textit{a priori} knowledge about the holes. Another modification
allows us to recover even the holes.

\section{Bergman polynomials on full domains vs. domains with holes}
\setcounter{equation}{0}
The following theorem shows that in many respect Bergman polynomials
on $G$ and on $G^\star$ behave similarly.
\begin{theorem}\label{th1}
If $G$ is a union of a finite family of bounded
Jordan domains lying a positive distance apart and $G^\star=G\setminus\mathcal{K}$, where $\mathcal{K}\subset G$ is compact,
 then, as $n\to\i$,
\begin{itemize}
\item[(a)]
\vspace{.2cm}
$\g_n(G^\star)/\g_n(G)\to 1$,
\item[(b)]
\vspace{.2cm}
$\|p_n(G^\star,\cdot)-p_n(G,\cdot)\|_{L^2(G)}\to 0,$
\item[(c)]
\vspace{.2cm}
$\l_n(G^\star,z)/\l_n(G,z)\to 1$ uniformly on compact subsets of $\ov\C\setm \ov G$,
\item[(d)]
\vspace{.2cm}
$p_n(G^\star,z)/p_n(G,z)\to 1$
uniformly on compact subsets of $\ov\C\setm {\rm Con}(G)$.
\end{itemize}
\end{theorem}

Here Con$(G)$ denotes the convex hull of $G$.

Since outside $\ov G$  both $\l_n(G^\star,z)$ and $\l_n(G,z)$ tend to zero
locally uniformly (see (\ref{9}) below), while inside $G$ both quantities
tend to a positive finite limit (see the next lemma), part (c) of Theorem \ref{th1}
is particularly useful in domain reconstruction (see Section 5), because
it tells us that, in the algorithm considered,
for reconstructing the outer boundary $\Gamma$ one does not need
to know in advance whether or not there are holes inside $G$.

The proof of Theorem \ref{th1} is based on
\begin{lem}\label{lem0} We have
\be \sum_{n=0}^\i |p_n(G^\star,z)|^2<\i\label{io}\ee
uniformly on compact subsets of $G$.
In particular, $p_n(G^\star,z)\to 0$ uniformly on compact subsets of $G$.
\end{lem}
\proofs. Let $V$ be a compact subset of $G$. Choose a system $\s\subset G^\star$ of
closed broken  lines
separating $V$ from $\partial G$ (meaning each $V\cap G_j$ is separated
from each $\partial G_j$),
 and choose $r>0$ such that
the disk $D_r(z)$ of radius $r$ about $z$ lies in $G^\star$ for all
$z\in \s$. For any $N>1$ and fixed $z\in \s$ we obtain from the subharmonicity in $t$ of
\[|P_N(t)|^2:=\left|\sum _{n=0}^N\overline{p_n(G^\star,z)}p_n(G^\star,t)\right|^2\]
the estimate
\bean \left(\sum_{n=0}^N|p_n(G^\star,z)|^2\right)^2&=&|P_N(z)|^2\le \frac{1}{r^2\pi}
\int_{D_r(z)}|P_N(t)|^2dA(t)\\
&\le&  \frac{1}{r^2\pi}\int_{G^\star}|P_N(t)|^2dA(t)
= \frac{1}{r^2\pi}\sum_{n=0}^N|p_n(G^\star,z)|^2.\eean
Thus,
\begin{equation}\label{eq:sumpN2}
\sum_{n=0}^N|p_n(G^\star,z)|^2\le  \frac{1}{r^2\pi}
\end{equation}
on $\s$, hence, again by  subharmonicity, the same
is true inside $\s$ (i.e. in every bounded component of
$\C\setm \s$). For  $N\to\infty$
we get
\be \sum_{n=0}^\i|p_n(G^\star,z)|^2\le  \frac{1}{r^2\pi}\label{fin}\ee
on and inside $\s$, but we still need to prove the
uniform convergence on $V$ of the series on the left hand side.

Let $\s_1$ be another family of closed broken lines
lying inside $\s$ separating $V$ and $\s$. If $\d$ is the
distance of $\s$ and $\s_1$, then for any $N$ and
any choice $|\e_n|=1$ we have, by Cauchy's formula for the derivative
of an analytic function for $z,w\in \s_1$
\bean
&&\left|\sum _{n=0}^N\e_n p_n(G^\star,z)p_n'(G^\star,w)\right|
\le \frac{L}{2\pi\d^2}\max_{t\in \s}\left|\sum _{n=0}^N\e_n p_n(G^\star,z)p_n(G^\star,t)\right|\\[10pt]
&&\quad \le  \frac{L}{2\pi\d^2}\max_{t\in \s}\left(\sum _{n=0}^N|p_n(G^\star,z)|^2\right)^{1/2}
\left(\sum _{n=0}^N|p_n(G^\star,t)|^2\right)^{1/2}
\le \frac{L}{2\d^2} \frac{1}{r^2\pi^2},
\eean
where $L$ is the length of $\s$.
So for $w=z$ an appropriate choice of the $\e_n$'s gives
\[
\sum _{n=0}^N |p_n(G^\star,z)||p_n'(G^\star,z)|\le \frac{L}{2\d^2} \frac{1}{r^2\pi^2}
\]
for all $z\in \s_1$. But then, if $ds$ is arc-length on $\s_1$, we obtain on $\s_1$
\[{\frac{d}{ds}\sum _{n=0}^N |p_n(G^\star,\cdot )|^2}\restr{\cdot=z}\le 2\sum_{n=0}^N|p_n(G^\star,z)||p_n'(G^\star,z)|
\le \frac{L}{\d^2} \frac{1}{r^2\pi^2},\]
which shows that on $\s_1$ the family
\[\left\{\sum _{n=0}^N |p_n(G^\star,z)|^2\right\}_{N=0}^\i\]
is uniformly equicontinuous. Since it converges pointwise
to a finite limit (see (\ref{fin})), we can conclude that the convergence in (\ref{fin})
is uniform on $\s_1$, and hence (by subharmonicity) also on $V$ (which lies inside
$\s_1$).\endproof

\proofs Theorem \ref{th1}.
In view of (\ref{eq:minimal2}) we have
\begin{equation}\label{9-}
\begin{alignedat}{1}
\frac{1}{\g_n(G)^2}&\le \int_G \frac{|p_n(G^\star,z)|^2}{\g_n(G^\star)^2}dA(z)=\int_{G^\star}+\int_\mathcal{K}\\
&\le \frac{1}{\g_n(G^\star)^2}+\frac{\e_n^2|\mathcal{K}|}{\g_n(G^\star)^2}=\frac{1+\e_n^2|\mathcal{K}|}{\g_n(G^\star)^2},
\end{alignedat}
\end{equation}
where
\be \label{09}
\e_n:=\|p_n(G^\star,\cdot)\|_\mathcal{K}\to 0%+\|p_n(G,\cdot)\|_\mathcal{K}\to 0
\ee
by Lemma \ref{lem0}. (Here and below we use $|\mathcal{K}|$ to denote the area measure of $\mathcal{K}$.)
On the other hand, (\ref{eq:comp-pri1}) gives that $\g_n(G^\star)\ge \g_n(G)$,
which, together with the preceding inequality shows
\begin{equation}\label{10}
1\le \frac{\g_n(G^\star)^2}{\g_n(G)^2}\le 1+\e_n^2|\mathcal{K}|,
\end{equation}
and this proves (a).
\bigskip

Next we apply a standard parallelogram-argument:
\bean \int_{G^\star}\left|\frac{1}{2}\left(\frac{p_n(G,\cdot)}{\g_n(G)}-\frac{p_n(G^\star,\cdot)}{\g_n(G^\star)}\right)\right|^2dA&+&
\int_{G^\star}\left|\frac{1}{2}\left(\frac{p_n(G,\cdot)}{\g_n(G)}+\frac{p_n(G^\star,\cdot)}{\g_n(G^\star)}\right)\right|^2dA\\[10pt]
&&\hskip-2cm=\frac{1}{2}\int_{G^\star}\left|\frac{p_n(G,\cdot)}{\g_n(G)}\right|^2dA+\frac{1}{2}\int_{G^\star}
\left|\frac{p_n(G^\star,\cdot)}{\g_n(G^\star)}\right|^2dA.\eean
By (\ref{eq:minimal2}) the second term on the left is $\ge 1/\g_n(G^\star)^2$, the second term on the right is
$1/(2\g_n(G^\star)^2)$ and, according to (\ref{9-}), the first term on the right is
\[\le
\frac12\int_{G}\left|\frac{p_n(G,\cdot)}{\g_n(G)}\right|^2dA=
\frac{1}{2\g_n(G)^2}\le \frac{1+\e_n^2|\mathcal{K}|}{2\g_n(G^\star)^2}.
\]
 Therefore, we can conclude
\[\int_{G^\star}\left|\frac{p_n(G,\cdot)}{\g_n(G)}-\frac{p_n(G^\star,\cdot)}{\g_n(G^\star)}\right|^2dA\le \frac{2\e_n^2|\mathcal{K}|}{\g_n(G^\star)^2},\]
and since (\ref{10}) implies
\[\left|1-\frac{\g_n(G^\star)}{\g_n(G)}\right|\le \e_n^2|\mathcal{K}|,\]
we arrive at
\begin{equation}\label{eqL2fromgamma}
\int_{G^\star}\left|p_n(G,\cdot)-p_n(G^\star,\cdot)\right|^2dA=O(\e_n^2),
\end{equation}
as $n\to\infty$. It is easy to see that the norms on $G^\star$ and $G$ for functions in $L^2_a(G)$ are equivalent; indeed, if $f\in L^2_a(G)$ and $\Gamma_0$ 
is the union of $m$ Jordan curves lying in $G^\star$ and containing $\mathcal{K}$ in its interior, then
$$ \|f\|_{L^2(G^\star)}^2 \leq \|f\|_{L^2(G)}^2=\|f\|_{L^2(G^\star)}^2+\|f\|_{L^2(\mathcal{K})}^2
$$
and, by subharmonicity, 
$$\|f\|_{L^2(\mathcal{K})}^2 \leq |\mathcal{K} |\max_{z \in \mathcal{K}}|f(z)|^2  \leq |\mathcal{K} |\max_{z \in \Gamma_0}|f(z)|^2\leq \frac{|\mathcal{K} |}{R^2\pi} \|f\|_{L^2(G^\star)}^2,
$$
where $R:=\dist(\Gamma_0,\partial G^\star ).$ Hence part (b) follows from (\ref{eqL2fromgamma}).

\bigskip

To prove (c), let $z$ lie in
$\ov \C\setm \overline G$. For an $\e>0$  select an $M$ such that
\be \sum_{j=M}^\i|p_j(G^\star,t)|^2\le \e, \qquad t\in \mathcal{K},\label{k1}\ee
(see Lemma \ref{lem0}). For the polynomial
\[P_n(t):=\frac{\sum _{j=M}^n\overline{p_j(G^\star,z)}p_j(G^\star,t)}{\sum _{j=M}^n|p_j(G^\star,z)|^2},\qquad n>M,\]
we have $P_n(z)=1$ and
\[\int_{G^\star}|P_n(t)|^2dA(t)=\frac{1}{\sum _{j=M}^n|p_j(G^\star,z)|^2}.\]
For its square integral over $\mathcal{K}$ we have by H\"older's inequality
\[\int_\mathcal{K} |P_n(t)|^2dA(t) \le\int_\mathcal{K} \frac{\sum_{j=M}^n|p_j(G^\star,t)|^2}{\sum _{j=M}^n|p_j(G^\star,z)|^2}dA(t)
\le \frac{|\mathcal{K}|\e}{\sum _{j=M}^n|p_j(G^\star,z)|^2}.\]
If we add together these last two integrals we obtain
\be \l_n(G,z)\le \frac{1+|\mathcal{K}|\e}{\sum _{j=M}^n|p_j(G^\star,z)|^2}.\label{oi}\ee

On the other hand, it is easy to see that outside $\ov G$ we always have
\be \sum _{j=0}^n|p_j(G^\star,z)|^2\to\i\label{9}\ee
as $n\to\i$, and actually this convergence to infinity is uniform on compact
subsets of $\Omega:=\ov\C\setm \ov G$. Indeed, if $\{F_n\}$ denotes a sequence of Fekete polynomials
associated with $\ov G$, then it is known (see e.g. \cite[Ch. III, Theorems 1.8, 1.9]{SaffTotik}) that
\begin{equation}\label{eq:Fek1}
\|F_n\|_{\ov G}^{1/n}\to \c(\ov G)=\c(\Gamma),\quad n\to\infty,
\end{equation}
where $\c(\ov G)$ denotes the \emph{logarithmic capacity} of $\ov G.$
At the same time
\begin{equation}\label{eq:Fek2}
|F_n(z)|^{1/n}\to \c(\ov G)\exp \left(g_{\Omega}(z,\i)\right),\quad n\to\infty,
\end{equation}
uniformly on compact subsets of  $\ov\C\setm \ov G$,
where $g_{\Omega}(z,\i)$ denotes the Green function of $\Omega$ with pole at infinity.
Thus,
\begin{equation}\label{eq:Fek3}
\l_n(G^\star,z)\le \int_{G^\star}\left|\frac{F_n(t)}{F_n(z)}\right|^2dA(t)\to 0,\quad n\to\infty,
\end{equation}
uniformly on compact subsets of  $\Omega$. (Note that
$g_{\Omega}(z,\i)$ has positive lower bound there.) Since $1/\l_n(G^\star,z)$
is the left-hand side of (\ref{9}), the relation (\ref{9}) follows.

Combining (\ref{oi}) and (\ref{9})
we can write
\bea \l_n(G^\star,z)\le \l_n(G,z)&\le& \frac{1+|\mathcal{K}|\e}{\sum _{j=M}^n|p_j(G^\star,z)|^2}
=(1+o(1))\frac{1+|\mathcal{K}|\e}{\sum _{j=0}^n|p_j(G^\star,z)|^2}\nonumber \\[10pt]
&=&(1+o(1))(1+|\mathcal{K}|\e)\l_n(G^\star,z),\label{k2}\eea
and since this relation is uniform on compact subsets of $\Omega$,
part (c) follows since $\e>0$ was arbitrary.
\bigskip

Finally, we prove part (d). Notice first of all that for $i,j\le n$
the expression $(z^it^j-z^jt^i)/(z-t)$ is a polynomial in $t$ of degree smaller than $n$,
therefore the same is true of
\[
\frac{p_n(G,z)p_n(G^\star,t)-p_n(G,t)p_n(G^\star,z)}{z-t},
\]
so this expression is orthogonal to $p_n(G,t)$ on $G$ with respect to area measure. Hence,
\[
\int_{G} \frac{p_n(G,z)p_n(G^\star,t)\ov{p_n(G,t)}}{z-t}dA(t)=
\int_{G}\frac{p_n(G,t)p_n(G^\star,z)\ov{p_n(G,t)}}{z-t}dA(t),
\]
and then division gives
\be
\frac{p_n(G^\star,z)}{p_n(G,z)}-1=
\frac{\int_{G}\frac{(p_n(G^\star,t)-p_n(G,t))\ov{p_n(G,t)}}{z-t}dA(t)}
{\int_{G}\frac{|p_n(G,t)|^2}{z-t}dA(t)}.\label{p1}
\ee

Let now $z$ be outside the convex hull of $G$ and let $z_0$ be
the closest point in the convex hull to $z$. Then $G$ lies
in the half-plane $\{t\sep \Re\{(z-t)/(z-z_0)\}\ge 1\}$, so for $t\in G$
\[\Re \frac{z-z_0}{z-t}= \frac{\Re\{(z-t)/(z-z_0)\}}{(|z-t|/|z-z_0|)^2}
\ge \frac{|z-z_0|^2}{|z-t|^2}\ge
\frac{|z-z_0|^2}{(|z-z_0|+{\rm diam}(G))^2}.\]
This gives the following bound for the modulus of the denominator in (\ref{p1}):
\bean
\left|{\int_{G}\frac{|p_n(G,t)|^2}{z-t}dA(t)}\right|&\ge& \frac{1}{|z-z_0|}\Re
{\int_{G}\frac{z-z_0}{z-t}|p_n(G,t)|^2dA(t)}\\[10pt]
&\ge& \frac{|z-z_0|}{(|z-z_0|+{\rm diam}(G))^2}
{\int_{G}|p_n(G,t)|^2dA(t)}\\[10pt]
&=&\frac{|z-z_0|}{(|z-z_0|+{\rm diam}(G))^2}.
\eean
On the other hand, in the numerator of (\ref{p1})
we have $1/|z-t|\le 1/|z-z_0|$, so we obtain from the
Cauchy-Schwarz inequality that
\bean
&&\left|\int_{G}\frac{(p_n(G^\star,t)-p_n(G,t))\ov{p_n(G,t)}}{z-t}dA(t)\right|\\
&&\hskip2cm\le \frac{1}{|z-z_0|}\left(\int_{G}|p_n(G^\star,t)-p_n(G,t)|^2dA(t)\right)^{1/2}.
\eean

Collecting these estimates we can see that
\[\left|\frac{p_n(G^\star,z)}{p_n(G,z)}-1\right|\le \frac{(|z-z_0|+{\rm diam}(G))^2}{|z-z_0|^2}
\|p_n(G^\star,\cdot)-p_n(G,\cdot)\|_{L^2(G)}.\]
Now invoking part (b), we can see that the left-hand side is uniformly small
on compact subsets of $\ov\C\setm {\rm Con}(G)$ since for ${\rm dist}(z,G)\ge \d$
we have
\[\frac{|z-z_0|+{\rm diam}(G)}{|z-z_0|}\le \frac{\d+{\rm diam}(G)}{\d}.\]
This proves (d)\footnote{The analysis used in the proof of part (d) was also found independently by B. Simanek (see ~\cite{Si12},
Lemma~2.1 and Theorem~2.2).}
\endproof

\section{Smooth outer boundary}
\setcounter{equation}{0}
Next, we make Theorem \ref{th1} more precise when the boundary $\G$ of $G$
is $C(p,\a)$-smooth, by which we mean that, for $j=1,\ldots,m,$  if $\g_j$ is the arc-length parametrization
of $\G_j$, then $\g_j$ is $p$-times differentiable, and its $p$-th derivative
belongs to the Lip $\a$.

Let $\|\cdot \|_{\overline G}$ denote the supremum norm on the closure $\overline G$ of $G$.
\begin{theorem}\label{th2}
If each of the boundary curves $\G_j$ is $C(p,\a)$-smooth for some $p\in \{1,2,\ldots\}$ and $0<\a<1$,
then
\begin{itemize}
\item[(a)]
$\g_n(G^\star)/\g_n( G)=1+O(n^{-2p+2-2\a}),$
\vspace{.2cm}
\item[(b)]
$\|p_n(G^\star,\cdot)-p_n(G,\cdot)\|_{\ov G}=O(n^{-p+2-\a}),$
\vspace{.2cm}
\item[(c)]
$\l_n(G^\star,z)/\l_n(G,z)=1+O(n^{-2p+3-2\a})$,  uniformly on compact subsets of $\ov\C\setm \ov G$,
\item[(d)]
\vspace{.2cm}
$p_n(G^\star,z)/p_n(G,z)=1+O(n^{-p+1-\a})$, uniformly on compact subsets of $\ov\C\setm {\rm Con}(G)$.
\end{itemize}
If each $\G_j$ is analytic, then {\textup{(a)}}--{\textup{(d)}} is true with
$O(q^n)$ on the right-hand sides for some $0<q<1$.
\end{theorem}
Note that now in (b) we have the supremum norm, so
$p_n(G^\star,z)-p_n(G,z)\to 0$ uniformly on $\ov G$ if $p>1$. Note also that nothing
like (d) is possible in the convex hull of $G$ since $p_n(G,\cdot)$ may have
zeros there, which need not be zeros of $p_n(G^\star,\cdot)$.

As background for the proof of Theorem 3.1, we shall first define $m$ special holes (lakes)
whose union contains $\mathcal{K}$. For this purpose, let
$\f_j$ map $G_j$ conformally onto the unit disk $\mathbb{D}$, and select an $0<r<1$ such
that each of the holes $\mathcal{K}_j:=\mathcal{K}\cap G_j$ is mapped by $\f_j$ into the disk $\mathbb{D}_r:=\{w: |w|<r\}$.
Let $\widetilde{\mathbb{D}}:=\{w: r<|w|<1\}$ and define $\widetilde G_j:=\f_j^{-1}(\widetilde{\mathbb{D}})$,
$\widetilde G:=\cup _{j=1}^m \widetilde G_j$.
Thus, the special holes  $\widetilde{\mathcal{K}}_j:=G_j\setm \widetilde G_j$ we are considering are
the preimages of the closed disk $\overline{\mathbb{D}}_r$ under $\f_j$.
Clearly, the above construction leads to the inclusions
\begin{equation}\label{eq:archi-inclusion}
\widetilde{G}\subset G^\star\subset G.
\end{equation}

We shall need to work with functions in  the Bergman space $L^2_a(G)$  but with
the inner product
\begin{equation}\label{eq:innpro-ti}
\langle f,g\rangle_{\widetilde{G}}:=\int_{\widetilde{G}}f(z) \overline{g(z)} dA(z),
\end{equation}
and corresponding norm $\|\cdot\|_{\widetilde{G}}.$
Let $L^{2\#}_a(G)$ denote the space of functions in $L^2_a(G)$ endowed  with the inner product (\ref{eq:innpro-ti}).
It is easy to see
that $L^{2\#}_a(G)$  is again a
Hilbert-space, but note that it is different from $L^2_a(\widetilde G)$ (the definition of the
norm on the two spaces is the same, but the latter space contains also functions
that may not be analytically continued throughout $G$, while the former space contains only
analytic functions in $G$). In fact, in $L^{2\#}_a(G)$, the polynomials
$\{p_n(\widetilde G,\cdot)\}_{n=0}^\i$ form a {\it complete} orthonormal system
(they also form an orthonormal system in $L^2_a(\widetilde G)$, which, however,
is not complete). Consequently, the reproducing kernel of $L^{2\#}_a(G)$ is
\begin{equation}\label{eq:Ksharp}
K^\#(z,\zeta)=\sum_{k=0}^\i \ov{p_k(\widetilde G,\zeta)}p_k(\widetilde G,z).
\end{equation}
Note that by Lemma \ref{lem0} (with $G^\star$ replaced by $\widetilde{G}$) the series on the right hand
side converges  uniformly on compact subsets of $G\times G$.

Analogously, we define the Hilbert space $L^{2\#}_a(\mathbb{D})$ consisting of functions in $L^{2}_a(\mathbb{D})$,
but with inner product
\begin{equation}%\label{eq:innpro-ti}
\langle f,g\rangle_{\widetilde{\mathbb{D}}}:=\int_{\widetilde{\mathbb{D}}}f(w) \overline{g(w)} dA(w).
\end{equation}

The following lemma provides a representation for the reproducing kernel $K^\#(z,\zeta)$ in terms of the reproducing
kernel for the space $L^{2\#}_a(\mathbb{D})$.
\begin{lemma}\label{lem:kercomp}
Let $J(w,\omega)$ denote the reproducing kernel for $L^{2\#}_a(\mathbb{D})$.
Then,
\begin{equation}\label{eq:KK1-ti}
K^{\#}(z,\zeta)=\left\{
\begin{array}{cl}
\overline{\varphi_j'(\zeta)}\varphi_j'(z)J(\varphi_j(z),\varphi_j(\zeta)), &\mathrm{if}\ z, \zeta\in{G}_j,\,\, j=1,\ldots,m, \\
 0,            &\mathrm{if}\ z\in{G}_j,\ \zeta\in{G}_k,\ j\neq k.
\end{array}
\right.
\end{equation}
Furthermore,
\begin{equation}\label{eq:Kw-ome}
J(w,\omega)=\sum_{\nu=0}^\i\frac{r^{2\nu}}{\pi (1-r^{2\nu}w\overline{\omega} )^2},\quad w,\omega\in\mathbb{D},
\end{equation}
and consequently, for $z,\zeta\in{G}_j$,
\begin{equation}\label{kfinal-ti}
K^{\#}(z,\zeta)=\overline{\varphi_j'(\zeta)}\varphi_j'(z)
\sum_{\nu=0}^\i\frac{r^{2\nu}}{\pi [1-r^{2\nu}\overline{\varphi_j(\zeta)}\varphi_j(z)]^2}.
\end{equation}
\end{lemma}
\proofs.
As with (\ref{eq:KK1}) it suffices to verify (\ref{eq:KK1-ti}) for $z,\zeta\in{G}_j$, $j=1,...,m$.
In fact, for $z,\zeta\in{G}_j$ the relation in (\ref{eq:KK1-ti}) is quite standard, see, e.g.,
\cite[Section 1.3, Theorem 3]{DuSc04}. To derive this relation, observe that
since the Jacobian of the mapping $w=\varphi_j(z)$  is $|\varphi_j'(z)|^2$, we have
\begin{equation*}
\int_{\widetilde G_j}|F(\varphi_j(z))|^2|\varphi_j'(z)|^2dA(z)=\int_{\widetilde{\mathbb{D}}} |F(w)|^2dA(w),
\end{equation*}
for any $F\in L^{2\#}_a(\mathbb{D})$.
Hence, the mapping $F\to F(\varphi_j)\varphi_j'$ is an isometry from $L^{2\#}_a(\mathbb{D})$ into
$L^{2\#}_a(G_j):=\{f{\chi}_{G_j}:f\in L^{2,\#}_a(G)\}$.
This mapping is actually onto $L^{2\#}_a(G_j)$, with inverse $f\to f(\f_j^{-1})(\f_j^{-1})'$.

Next, from the reproducing property of $J(w,\omega)$, it follows that for $\omega\in \mathbb{D}$,
\begin{equation*}
F(\omega)=\int_{\widetilde{\mathbb{D}}} F(w)\overline{J(w,\omega)}dA(w),\quad
F\in L^{2\#}_a(\mathbb{D}).
\end{equation*}
If we make the change of variable $w=\varphi_j(z)$, $\omega=\varphi_j(\zeta)$, this takes the form
\begin{equation*}
F(\varphi_j(\zeta))=\int_{\widetilde G_j} F(\varphi_j(z))\overline{J(\varphi_j(z),\varphi_j(\zeta))}|\varphi_j'(z)|^2dA(z),
\quad \zeta\in G_j,
\end{equation*}
which, after multiplication by $\varphi_j'(\zeta)$ gives for $f(\zeta):=F(\varphi_j(\zeta))\varphi_j'(\zeta)$ that
\begin{equation}\label{eq:KK-Gj}
f(\zeta)=\int_{\widetilde G_j}
f(z)\overline{\overline{\varphi_j'(\zeta)}\varphi_j'(z)J(\varphi_j(z),\varphi_j(\zeta))}dA(z),
\quad \zeta\in G_j.
\end{equation}
Thus $\overline{\varphi_j'(\zeta)}\varphi_j'(z)J(\varphi_j(z),\varphi_j(\zeta))$ is the reproducing kernel for the space
$L^{2\#}_a(G_j)$, which establishes (\ref{eq:KK1-ti}).

To obtain the formula for $J(w,\omega)$, we note that the polynomials
\begin{equation*}
\left(\frac{\pi}{n+1}\left(1-r^{2n+2}\right)\right)^{-1/2}w^n,\quad n=0,1,\ldots,
\end{equation*}
form a complete orthonormal system in the space $L^{2\#}_a(\mathbb{D})$.
Therefore, we obtain the following representation:
\begin{equation*}%\label{}
\begin{alignedat}{1}
J(w,\omega)&=\sum_{n=0}^\i \left(\frac{\pi}{n+1}\left(1-r^{2n+2}\right)\right)^{-1}w^n\overline{\omega}^n
=\sum_{n=0}^\i \frac{n+1}{\pi}\sum_{\nu=0}^\i r^{2\nu}r^{2n\nu}w^n\overline{\omega}^n\\
&=\sum_{\nu=0}^\i r^{2\nu}\sum_{n=0}^\i \frac{n+1}{\pi} r^{2n\nu}w^n\overline{\omega}^n
=\sum_{\nu=0}^\i\frac{r^{2\nu}}{\pi (1-r^{2\nu}w\overline{\omega} )^2},
\end{alignedat}
\end{equation*}
and the result (\ref{kfinal-ti}) follows from (\ref{eq:KK1-ti}).
\endproof

\noindent
\textbf{Proof of Theorem \ref{th2}}.
With the above preparations we now turn to the proof of part (a) in Theorem \ref{th2}.
First, we need a good polynomial
approximation of the kernel $K^{\#}(\cdot,\zeta)$ on $\overline{G}$, for fixed $\zeta\in V$, where
$V$ {is a compact subset  of} $G_j$.
By the Kellogg-Warschawskii theorem (see, e.g., \cite[Theorem 3.6]{Pommerenke}),
our assumption $\Gamma_j\in C(p,\alpha)$ implies that $\varphi_j$ belongs to the class $C^{p+\alpha}$ on $\Gamma_j$.
Thus, $\varphi'_j\in C^{p-1+\alpha}$ on $\Gamma_j$ and (\ref{kfinal-ti}) shows that the kernel
$K^{\#}(\cdot,\zeta)$ is a $C^{p-1+\alpha}$-smooth function on $\Gamma_j$
and the smoothness is uniform when $\zeta$ lies in a compact subset $V$ of $G_j$.
Consequently  (see, e.g., \cite[p.~34]{Su74}),  there are
polynomials $P_{\nu,j,\zeta}(z)$ of degree $\nu$ such that for $\zeta\in V$
\begin{equation*}%\label{knt0}
\sup_{z\in \Gamma_j}|K^{\#}(z,\zeta)-P_{\nu,j,\zeta}(z)|\le C(\Gamma_j,V)\frac{1}{\nu^{p-1+\alpha}},
\quad \nu\in\mathbb{N},\,j=1, \ldots,m,
\end{equation*}
where $C(\Gamma_j,V)$ here and below denotes a positive constant, not necessarily the
same at each appearance,  that depends on $\Gamma_j$ and $V$, but is independent of $\nu$. Therefore, the maximum modulus principle gives
\begin{equation}\label{knt}
\sup_{z\in\overline{G}_j}|K^{\#}(z,\zeta)-P_{\nu,j,\zeta}(z)|\le C(\Gamma_j,V)\frac{1}{\nu^{p-1+\alpha}},
\quad \zeta\in V.
\end{equation}

Note that this provides a good approximation to $K^{\#}(z,\zeta)$ only for $z\in\overline{G}_j$.
However, $K^{\#}(z,\zeta)$ is also defined for $z\in\overline{G}_k$, $k\not=j$.
Actually, as we have seen in (\ref{eq:KK1-ti}), for such values $K^{\#}(z,\zeta)=0$. Therefore, in order to
obtain a good approximation to $K^{\#}(z,\zeta)$ for all $z\in\overline{G}$,
we have to modify the polynomials $\{P_{\nu,j,\zeta}(z)\}$. To this end, we note that since (\ref{knt})
implies that the $\{P_{\nu,j,\zeta}(z)\}$ are bounded uniformly for $z\in\overline{G}_j$, $\zeta\in V$ and $\nu\ge 1$,
the Bernstein-Walsh lemma~\cite[p.~77]{Walsh} implies that there is a constant $\tau>0$ such that
\begin{equation}\label{knt1}
|P_{\nu,j,\zeta}(z)|\le C(\Gamma,V)\tau^\nu,\quad z\in\overline{G}.
\end{equation}

Consider next the characteristic function
\begin{equation}\label{eq:chi-def}
\chi_{\overline{G}_j}(z):=\left\{
\begin{array}{cl}
1, &\mathrm{if}\ z\in\overline{G}_j,\\
0, &\mathrm{if}\ z\in\overline{G}_k,\  k\neq j.
\end{array}
\right.
\end{equation}
Since $\chi_{\overline{G}_j}$ has an analytic continuation to an open set containing $\overline{G}$, it is known from the theory of polynomial
approximation (cf. \cite[p. 75]{Walsh})
that there exist polynomials $H_{n/2,j}(z)$ of degree at
most $n/2$ such that \begin{equation}\label{knt2}
\sup_{z\in\overline{G}}|\chi_{\overline{G}_j}(z)-H_{n/2,j}(z)|\le  C(\Gamma,V)\eta^n,
\end{equation}
for some $0<\eta<1$.

For some small $\epsilon>0$ we set
$$
Q_{n,j,\zeta}(z):=P_{\epsilon n,j,\zeta}(z)H_{n/2,j}(z).
$$
This is a polynomial in $z$ of degree at most $\epsilon n +(n/2)<n$,
and (\ref{eq:chi-def})--(\ref{knt2}), in conjunction with (\ref{knt})--(\ref{knt1}), yield for large $n$
\begin{equation*}
\sup_{z\in\overline{G}_j}|K^{\#}(z,\zeta)-Q_{n,j,\zeta}(z)|\le C(\Gamma_j,V)\frac{1}{(\epsilon n)^{p-1+\alpha}}
+ C(\Gamma,V)\tau^{\epsilon n}\eta^n,
\end{equation*}
and
\begin{equation*}
\sup_{z\in\overline{G}\setminus\overline{G}_j}|K^{\#}(z,\zeta)-Q_{n,j,\zeta}(z)|\le  C(\Gamma,V)\tau^{\epsilon n}\eta^n,
\quad \zeta\in V\subset G_j.
\end{equation*}
Thus, if we fix $\epsilon>0$ so small that $\tau^\epsilon\eta<1$ is satisfied,
 we obtain for large enough $n$
\begin{equation}\label{KtiG-Qn}
\sup_{z\in\overline{G}}|K^{\#}(z,\zeta)-Q_{n,j,\zeta}(z)|\le C(\Gamma,V)\frac{1}{n^{p-1+\alpha}}.
\end{equation}
This is our desired estimate.

Since $Q_{n,j,\zeta}(z)$ is of degree smaller than $n$, using the reproducing property of the
kernel $K^{\#}(z,\zeta)$ and the orthonormality of $p_n(\widetilde{G},z)$ with respect to the inner product (\ref{eq:innpro-ti}),
we conclude that
\begin{equation*}%\label{}
\begin{alignedat}{1}
p_n(\widetilde G,\zeta)&=\langle p_n(\widetilde{G},\cdot), K^{\#}(\cdot,\zeta)\rangle_{\widetilde{G}}\\
&=\langle p_n(\widetilde{G},\cdot), K^{\#}(\cdot,\zeta)-Q_{n,j,\zeta}\rangle_{\widetilde{G}}.
\end{alignedat}
\end{equation*}
Therefore, from the Cauchy-Schwarz inequality and (\ref{KtiG-Qn}), we obtain the following uniform estimate for
$\zeta\in V$:
\begin{equation*}
|p_n(\widetilde G,\zeta)|\le C(\Gamma,V)\frac{1}{n^{p-1+\alpha}},
\end{equation*}
where we recall that $V$ is a compact subset of $G_j$.
Since this is true for any $j=1,\dots,m$, we have  shown that
\begin{equation}\label{k3}
|p_n(\widetilde G,\zeta)|\le C(\Gamma,V)\frac{1}{n^{p-1+\alpha}},\quad\zeta\in V,
\end{equation}
where now $V$ is any compact subset of $G$.

Consequently, with $V=\widetilde{\mathcal{K}}:=\cup_{j=1}^m\widetilde{\mathcal{K}}_j$ in (\ref{k3}),
and $G^\star$ and $\mathcal{K}$ replaced by $\widetilde{G}$ and $\widetilde{\mathcal{K}}$ in (\ref{9-}) and
(\ref{09}), from (\ref{10}) we get
\begin{equation}\label{eq:gamm-ti-gamma}
\frac{\gamma_n(\widetilde{G})}{\gamma_n(G)}=1+O\left(\frac{1}{n^{2(p-1+\alpha)}}\right),
\end{equation}
which in view of  the fact
\begin{equation*}%\label{eq:comp-pri-gamma-all}
\gamma_n(G)\le\gamma_n(G^\star)\le\gamma_n(\widetilde{G}),
\end{equation*}
yields part (a) of the theorem.

To prove part (b), notice that (\ref{eq:gamm-ti-gamma}) is (\ref{10}) with $\varepsilon_n=O(n^{-p+1-\alpha})$,
and so the argument leading from (\ref{10}) to (\ref{eqL2fromgamma}) yields
\begin{equation}\label{eq:L2Gs-esti}
\|p_n(G,\cdot)-p_n(G^\star,\cdot)\|_{L^2(G^\star)}=O\left(\frac{1}{n^{p-1+\alpha}}\right).
\end{equation}
The $L^2$-estimate in (\ref{eq:L2Gs-esti}) holds also over $G$ since, as was previously remarked,  the two norms
$\|\cdot\|_{L^2(G)}$ and $\|\cdot\|_{L^2(G^\star)}$ are equivalent in $L^2_a(G)$. %cf. (\ref{eq:norm-equiv}).
The uniform estimate in part (b) then follows from the $L^2$-estimate by using the  inequality
$$
\|Q_n\|_{\ov G}\le C(\Gamma)n\|Q_n  \|_{L^2(G)},
$$
which is valid for all polynomials $Q_n$ of degree at most $n\in\mathbb{N}$, where the constant $C(\Gamma)$
depends on $\Gamma$  only; see \cite[p. 38]{Su74}.

In proving part (c)  we may assume $p+\a> 3/2$ (see Theorem \ref{th1} (c)).
It follows from (\ref{k3})
that
\[\sum_{k=n}^\i |p_k(\widetilde G,z)|^2=O(n^{-2p+3-2\a})\]
uniformly on compact subsets of $G$, i.e. (\ref{k1}) holds (for $\widetilde G$ in place
of $G^\star$) with
$\e=O(n^{-2p+3-2\a})$. Copying the proof leading from (\ref{k1}) to (\ref{k2})
with this $\e$ we get
\[\l_n(\widetilde G,z)\le \l_n(G,z)=(1+O(n^{-2p+3-2\a}))\l_n(\widetilde G,z)\]
(indeed, by that proof the $o(1)$ in (\ref{k2}) is exponentially small).
In view of $\widetilde G\subset G^\star\subset G$ this then implies
\bean \l_n(G^\star,z)&\le& \l_n(G,z)=(1+O(n^{-2p+3-2\a}))\l_n(\widetilde G,z)\\[10pt]
&\le& (1+O(n^{-2p+3-2\a}))\l_n(G^\star,z),\eean
which is part (c) in the theorem.

Part (d) follows at once from the $L^2$-estimate in (\ref{eq:L2Gs-esti}), by working as in the proof of (d) in Theorem~\ref{th1}.

Regarding the case when all the curves $\Gamma_j$ are analytic, we have that
the conformal maps  $\f_j$ are analytic on $\overline{G}_j$,
and then so is the kernel $K^{\#}(z,\zeta)$ for $z\in \ov G$, and all fixed $\zeta\in\widetilde G$.
More precisely, if $V$ is a compact subset of $\widetilde G$, then
there is an open set $\ov G\subset U$ such that for $\zeta\in V$ the kernel $K(z,\zeta)$
is analytic for $z\in U$.
Then, from the proof of the classical polynomial approximation theorem for analytic functions mentioned previously,
together with the formula for $K^{\#}(z,\zeta)$, it follows that
there is a $0<q<1$ and a constant $C$ independent of $\zeta\in V$, such that in place of (\ref{knt}) we have
\be \sup_{z\in \tilde G_j}|K^{\#}(z,\zeta)-P_{n-1,j,\zeta}(z)|\le Cq^n, \quad \zeta\in V.
\label{knt*}\ee
Thus, instead of (\ref{k3}), we obtain
\bean
|\ov{p_n(\widetilde G,\zeta)}|&=&\left|\int_{\widetilde G} K^{\#}(z,\zeta)\ov{p_n(\widetilde G,z)}\,dA(z)\right|
\\&=&
\left|\int_{\widetilde G} (K^{\#}(z,\zeta)-P_{n-1,j,\zeta}(z))\ov{p_n(\widetilde G,z)}\,dA(z)\right|
\le  C|G|^{1/2}q^n,
\eean
so the $\e_n$ in (\ref{09}) is $O(q^n)$, and then the proofs of (a)--(d) above give
the same statements with error $O(q^n)$ (for a possibly different $0<q<1$).
\endproof

%\medskip
\begin{remark}
\textup{Our theorems thus far have emphasized the similar asymptotic behavior of the Bergman orthogonal polynomials
for an archipelago without lakes and the Bergman polynomials for an archipelago with lakes. Differences appear,
however, when one considers the asymptotic behaviors of the zeros of the two sequences of polynomials.
A future paper will be devoted to this topic.}
\end{remark}

\section{Asymptotics behavior}\label{section:asym-beha}
\setcounter{equation}{0}

Since area measure on the archipelago $G$ belongs to the class \textbf{Reg} of measures (cf. \cite{StahlTotik}),
it readily follows from Theorem~\ref{th1} that so does area measure on $G^\star$. In particular,
\begin{equation}\label{eq:nthlead*}
\lim_{n\to\infty}\gamma_n(G^\star)^{1/n}=\frac{1}{\capGm}.
\end{equation}

In order to describe the $n$-th root asymptotic behavior for the Bergman polynomials
$p_n(G^\star,z)$ in $\Omega$, we need the Green function $g_\Omega(z,\infty)$ of $\Omega$ with pole at infinity.
We recall that $g_\Omega(z,\infty)$ is harmonic
in $\Omega\setminus\{\infty\}$, vanishes on the boundary $\Gamma$ of $G$ and near $\infty$ satisfies
\begin{equation}\label{green_inf}
g_\Omega(z,\infty)=\log|z|+\log\frac{1}{\capGm}+O\left(\frac{1}{|z|}\right),
\quad |z|\to\infty,
\end{equation}
Our next result corresponds to Proposition 4.1 of \cite{GPSS} and follows in a similar manner.

\begin{proposition}\label{pro:nthroot}
%[\cite{GPSS}, Proposition 3.1]
The following assertions hold:
\begin{itemize}
\item[(a)] For every $z\in \ov{\C}\setminus{\rm{Con}}(G)$
and for any $z\in{\rm{Con}}(G)\setminus\Gbar$ not a limit point of zeros of the $p_n(G^\star,\cdot)$'s, we have
\begin{equation}\label{eq:nthroot3}
\lim_{n \to\infty}|p_n(G^\star,z)|^{1/n}=\exp\{g_\Omega(z,\infty)\}.
\end{equation}
The convergence is uniform on compact subsets of
$\ov{\C}\setminus{\rm{Con}}(G)$.
\item[(b)] There holds
\begin{equation}\label{eq:nthroot2}
\limsup_{n \to\infty}|p_n(G^\star,z)|^{1/n}=\exp\{g_\Omega(z,\infty)\},\ z\in\overline{\Omega},
\end{equation}
locally uniformly in ${\Omega}$.
\end{itemize}
\end{proposition}

For our next result we assume that all the boundary curves $\Gamma_j$ are analytic.
Its proof is a simple consequence of Theorem~4.1 of \cite{GPSS} in conjunction with  Theorem~\ref{th2} above.
\begin{proposition}\label{th:GPSS_ln}
Assume that every curve $\Gamma_j$, $j=1,\ldots,m$, constituting $\Gamma$ is analytic.
Then there exist positive constants $C_1(\Gamma,\mathcal{K})$ and $C_2(\Gamma,\mathcal{K})$ such that
\begin{equation}\label{eq:GPSS_lnge}
C_1(\Gamma,\mathcal{K})\le\sqrt{\frac{n+1}{\pi}}
\frac{1}{\gamma_n(G^\star)\,\capGm^{n+1}}\le C_2(\Gamma,\mathcal{K}),\quad n\in\mathbb{N}.
\end{equation}
\end{proposition}

As the following example emphasizes, we cannot expect that the limit of the sequence in (\ref{eq:GPSS_lnge})
exists when $m\ge 2$.
\begin{example}[\cite{GPSS}, Remark 7.1]
Consider the $m$-component lemniscate $G:=\{z:|z^m-1|<r^m\}$, $m\ge 2$, $0<r<1$, for which ${\textup{cap}(\Gamma)}=r$.
Then, the sequence
$$
\sqrt{\frac{n+1}{\pi}}
\frac{1}{\gamma_n(G)\,\capGm^{n+1}},\quad n\in\mathbb{N},
$$
has exactly m limit points:
$$
{r^{m-1}}, {r^{m-2}}, \ldots, {r}, 1.
$$
\end{example}

Combining the result of Theorem~\ref{th2} with that of Theorem 4.4 of \cite{GPSS}, we arrive at
estimates for the Bergman polynomials $\{p(G^\star,z\}$ in the exterior domain $\Omega$, where
we use $\dist(z,E)$ to denote the (Euclidean) distance of $z$ from a set $E$.
\begin{theorem}\label{th:PnOme}
With $G$ as in Proposition~\ref{th:GPSS_ln}, the following hold:
\begin{itemize}
\item[(a)]
There exists a positive constant $C$, such that
\begin{equation}\label{eq:pn-below}
|p_n(G^\star,z)|\le\ \frac{C}{\textup{dist}(z,\Gamma)}\sqrt{n}\exp\{ng_\Omega(z,\infty)\},\ \ z
\notin \overline{G}.
\end{equation}
\item[(b)]
For every $\varepsilon>0$ there exist a constant $C_{\varepsilon}>0$, such that
\begin{equation}\label{eq:pn-above}
|p_n(G^\star,z)|\ge\ C_{\varepsilon} \sqrt{n}\exp\{ng_\Omega(z,\infty)\},\quad
\textup{dist}(z,{\rm Con}(G))\ge\varepsilon.
\end{equation}
\end{itemize}
\end{theorem}

\section{Reconstruction algorithm from moments}\label{section:reconsrt}
\setcounter{equation}{0}
The present section contains the description and analysis of a reconstruction algorithm for the
archipelago with lakes $G^\star$, for the case when the lakes are themselves finite unions of disjoint Jordan regions.
The algorithm is motivated by the `reconstruction from moments' algorithm of \cite[Section 5]{GPSS} and
the estimates established in the previous sections.
In \cite{GPSS} the functional $\lambda^{1/2}_n(G,z)$ was used as the main reconstruction tool for recovering the shape of the
archipelago $G$ using area complex moment measurements.
Here we describe how to recover from $\lambda^{1/2}_n(G^\star,z)$ both the shape of $G$ and of its lakes.

Assume that the following set of area complex moments is available:
$$
\mu^\star_{ij}:=\int_{G^\star}
z^i\overline{z}^j\,dA(z),\quad i,j=0,1,\ldots,n.
$$
(For a discussion of how these moments are related to the real moments
$$\tau^\star_{m n}:=\int_{G^\star}x^my^n dxdy
$$
that arise in geometric tomography from measurements of the Radon transform, see \cite{GPSS} and \cite{SafSty}.) Our algorithm consists of two phases.

\medskip
\noindent {\tt RECONSTRUCTION ALGORITHM}
\newline\noindent
{\tt Phase A: Recovery of $G$}
\begin{enumerate}[I]
\itemsep=5pt
\item
Use the Arnoldi Gram-Schmidt process described below to compute $p_0(G^\star,z)$,
$p_1(G^\star,z),\ldots,p_n(G^\star,z)$,
from the given set of moments $\mu^\star_{i,j}$ of $G^\star$, $i,j=0,1,\ldots,n$.
\item
Plot the zeros of $p_n(G^\star,z)$.
\item
Form $\lambda^{1/2}_n(G^\star,z)$.
\item
Plot the level curves of the function $\lambda^{1/2}_n(G^\star,x+iy)$ on a
suitable rectangular frame for $(x,y)$ that surrounds the plotted
zero set\footnote{See Remark~\ref{rem:fejer}.}.
The outer-most level curves will provide an approximation to the boundary of $G$.
Denote by $\widehat{G}$ the region(s) bounded by this approximation.
\end{enumerate}
\noindent
{\tt Phase B: Recovery of $\mathcal{K}$}
\begin{enumerate}[I]
\itemsep=5pt
\item
Use the approximation  $\widehat{G}$ of $G$ to calculate the  moments
$$
\widehat{\mu}_{i,j}:=\int_{\widehat{G}} z^i\overline{z}^j\,dA(z), \quad i,j=0,1,\ldots,n.
$$
\item
Compute the approximate moments $\mu^\prime_{i,j}$ for the lakes $\mathcal{K}$ by taking the difference
$\widehat{\mu}_{i,j}-\mu^\star_{i,j}$
\item
Repeat steps I-IV of {\tt Phase A} with data $\mu^\prime_{i,j}$ in the place of $\mu^\star_{i,j}$,
to produce an approximation $\widehat{\mathcal{K}}$ to $\mathcal{K}$.
\end{enumerate}

\medskip
Step \textrm{I} of \texttt{Phase B} is computationally demanding, but can be carried out by approximating
the outer-most level curves by polygonal curves which will facilitate the computation of the area moments of $\widehat{G}$.
This aspect of the algorithm will be explored in a future paper. Here, we shall illustrate our method by using the
moments of $G$ instead of $\widehat{G}$.

We recall that the Gram-Schmidt (GS) process (mentioned in step I) converts, in an iterative fashion, a set of linearly
independent functions in some inner product space into a set of orthonormal polynomials
$\{p_0,p_1,\ldots,p_{n-1},p_{n}\}$.
By the \textit{Arnoldi} GS we mean the application of the GS process in the following way:
At the $k$-step, where the orthonormal polynomial $p_k$ is to be constructed, we use the polynomials
$\{p_0,p_1,\ldots,p_{k-1},zp_{k-1}\}$ as input of the process.
We refer to \cite[Section 7.4]{St-CA13} for a discussion regarding the stability properties of the Arnoldi GS.
In particular, we note that the Arnoldi  GS does not suffer from the severe
ill-conditioning associated with the conventional GS as reported, for instance, by theoretical and numerical
evidence in \cite{PW86}.

\begin{remark}\label{rem:fejer}
\textup{A well-known result of Fej\'{e}r asserts that the zeros of orthogonal polynomials with respect to a compactly supported
measure are contained in the convex hull of the support of the measure.
Thus the frames chosen in \texttt{Phases A} and \texttt{B} should at least contain such zeros. However, adjustments to the size of such
frames may be required, as may be indicated by the appearance of level lines for $\lambda_n^{1/2}$ that are not closed
(see Figure \ref{fig:rec_disk-pent-disk80b-e}).}
\end{remark}

The following theorem contains estimates for the asymptotic behavior of $\lambda^{1/2}_n(G^\star,z)$,
thus providing the theoretical support of the reconstruction algorithm given above.

\begin{theorem}\label{th:lambdan-Ome}
Under the general assumption that $\Gamma$ consists of a finite union of Jordan curves we have the following:
\begin{itemize}
\item[(a)]
There exists a positive constant $C$ such that
\begin{equation}\label{eq:th:lambdan-Ome-1}
\lambda^{1/2}_n(G^\star,z)\ge C\dist(z,\Gamma),\quad z\in G.
\end{equation}
\item[(b)]
For every compact subset $B$ of $\Omega$, there exists a positive constant $C(B)$ such that
\begin{equation}\label{eq:th:lambdan-Ome-2}
\lambda^{1/2}_n(G^\star,z)\le C(B)\exp\{-ng_\Omega(z,\infty)\},\quad z\in B.
\end{equation}
\end{itemize}
\end{theorem}

The estimate in (\ref{eq:th:lambdan-Ome-1}) is immediate from (\ref{eq:sumpN2}),
while (\ref{eq:th:lambdan-Ome-2}) follows from (\ref{eq:Fek1}) and (\ref{eq:Fek2}).

Regarding the use of the square root $\lambda_n^{1/2}$ rather than $\lambda_n$ itself,
as indicated in (\ref{eq:th:lambdan-Ome-1}), the former quantity decays linearly to zero with
the distance to the boundary $\Gamma=\partial G$, while the latter has a more rapid decay which will effect the omission
(due to negligibility) of level curves that are closer to $\Gamma$.
This can be seen by comparing Figure~\ref{fig:diskpen-vili} with the more accurate Figure~\ref{fig:diskpenA},
where the Maple routine \texttt{contourplot} was used to generate the level curves.

\begin{figure}[h]
\begin{center}
\includegraphics*[scale=0.4]{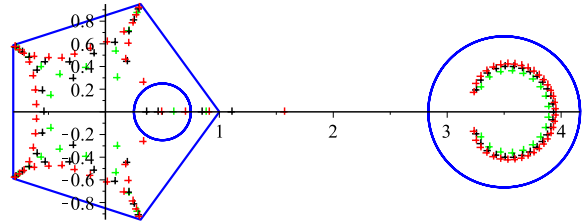}
\end{center}
\caption{Zeros of the polynomials $p_n(G^\star,z)$ of Example~\ref{ex:diskpen}, for $n=40,60$ and $80$.}
\label{fig:ro-diskpent80}
\end{figure}

\begin{figure}[h]
\begin{center}
\includegraphics*[scale=0.4]{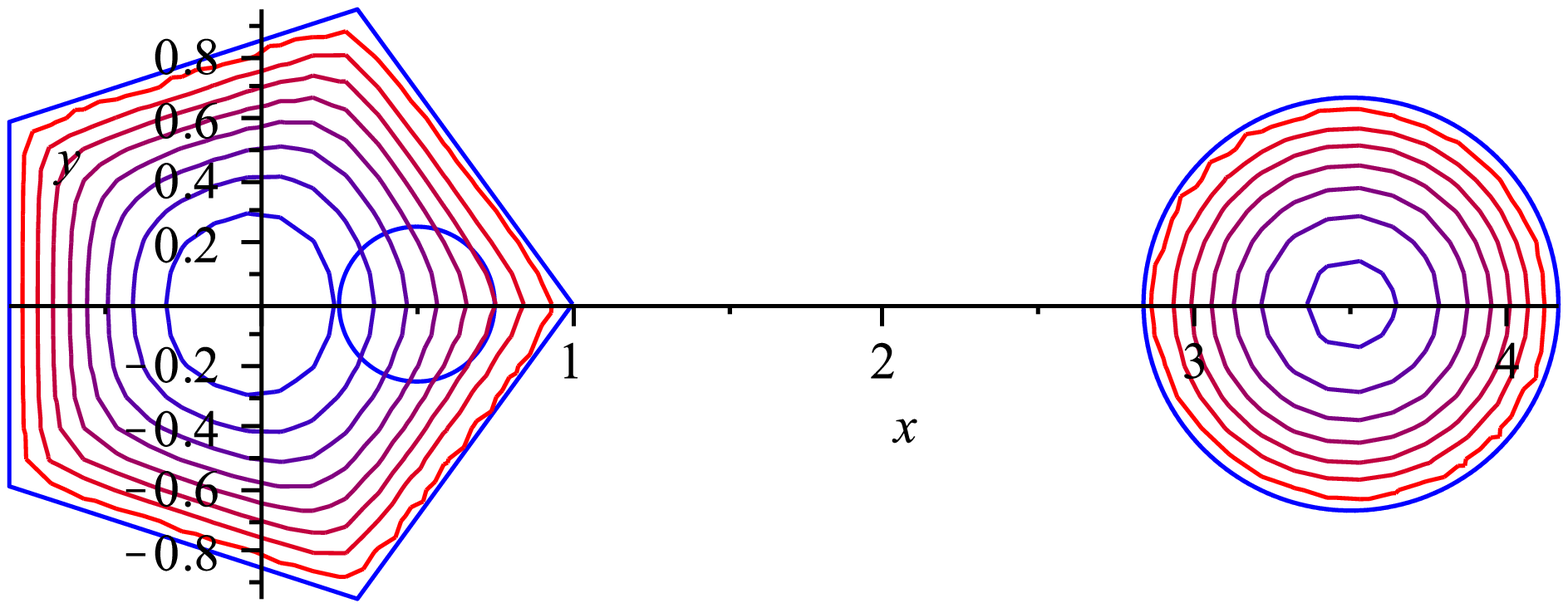}
\end{center}
\caption{\texttt{Phase A:} Level curves of $\lambda^{1/2}_{80}(G^\star,x+iy)$, on $\{(x,y):-2\le
x\le 5,-2\le y\le 2\}$, with $G^\star$ as in Example~\ref{ex:diskpen}.}
\label{fig:diskpenA}
\end{figure}

\begin{example}\label{ex:diskpen}
Recovery for the archipelago $G=G_1\cup G_2$, with $G_1$
denoting the canonical pentagon with vertices at the fifth roots of unity, $G_2=\{z:|z-7/2|<2/3\}$,
and lake $\mathcal{K}$ the closed disc centered at $1/2$ with radius $1/4$.
The boundaries of the archipelago $G^\star:=G\setminus\mathcal{K}$ are depicted in Figure \ref{fig:ro-diskpent80}.
\end{example}

In view of Remark~\ref{rem:fejer}, the zeros of the polynomial $p_n(G^\star,z)$ will give an indication of
the position of $G$ in the complex plane.
Accordingly, in Figure~\ref{fig:ro-diskpent80} we show the zeros for $n=40,60$ and $80$.
This should be compared with Figure 8 in \cite{GPSS}, which depicts zeros of $p_n(G,z)$.

\begin{figure}[h]
\begin{center}
\includegraphics*[scale=0.4]{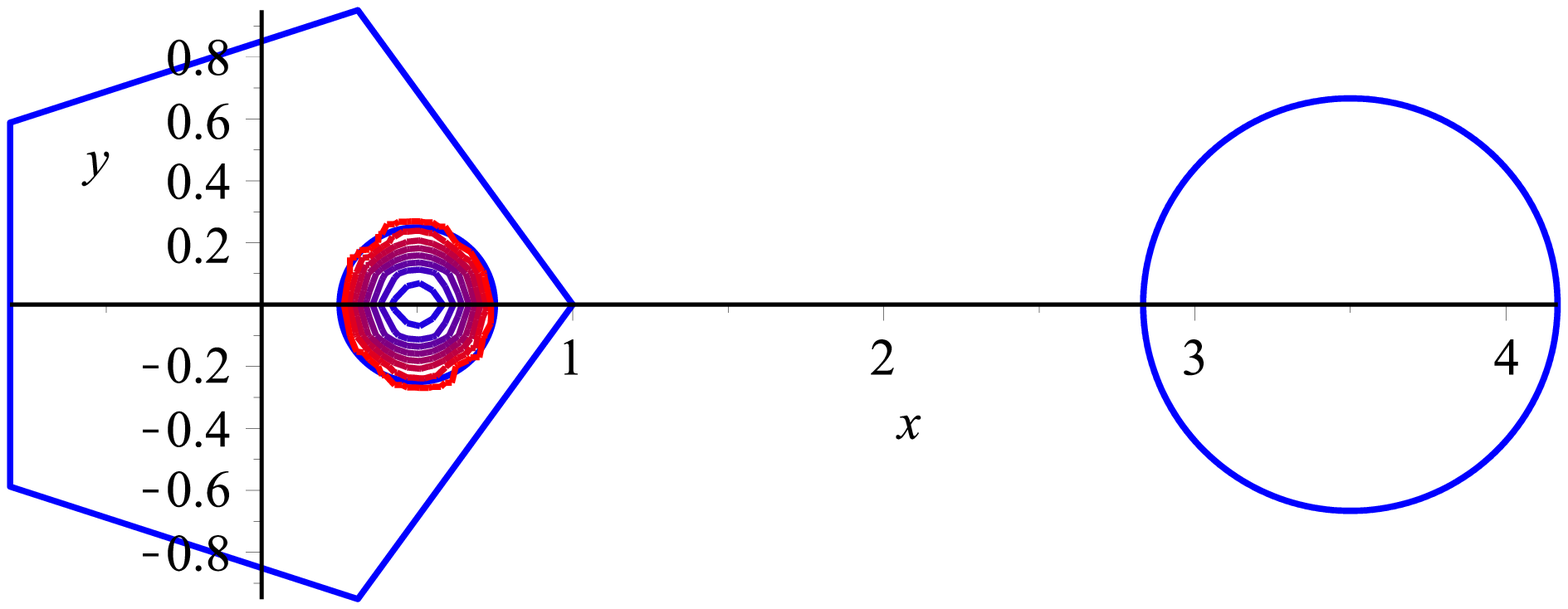}
\end{center}
\caption{\texttt{Phase B:} Level curves of $\lambda^{1/2}_{80}(\widehat{\mathcal{K}},x+iy)$, on $\{(x,y):-2\le
x\le 5,-2\le y\le 2\}$, with $G^\star$ as in Example~\ref{ex:diskpen}.}
\label{fig:diskpenB}
\end{figure}

\begin{figure}[h]
\begin{center}
\includegraphics*[scale=0.4]{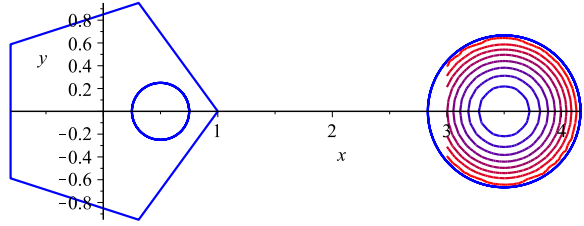}
\end{center}
\caption{\texttt{Phase A:} Level curves of $\lambda^{1/2}_{80}(G^\star,x+iy)$, for the inappropriately frame
 $\{(x,y):3\le
x\le 6,-2\le y\le 2\}$, with $G^\star$ as in Example~\ref{ex:diskpen}.}
\label{fig:rec_disk-pent-disk80b-e}
\end{figure}

\begin{figure}[h]
\begin{center}
\includegraphics*[scale=0.4]{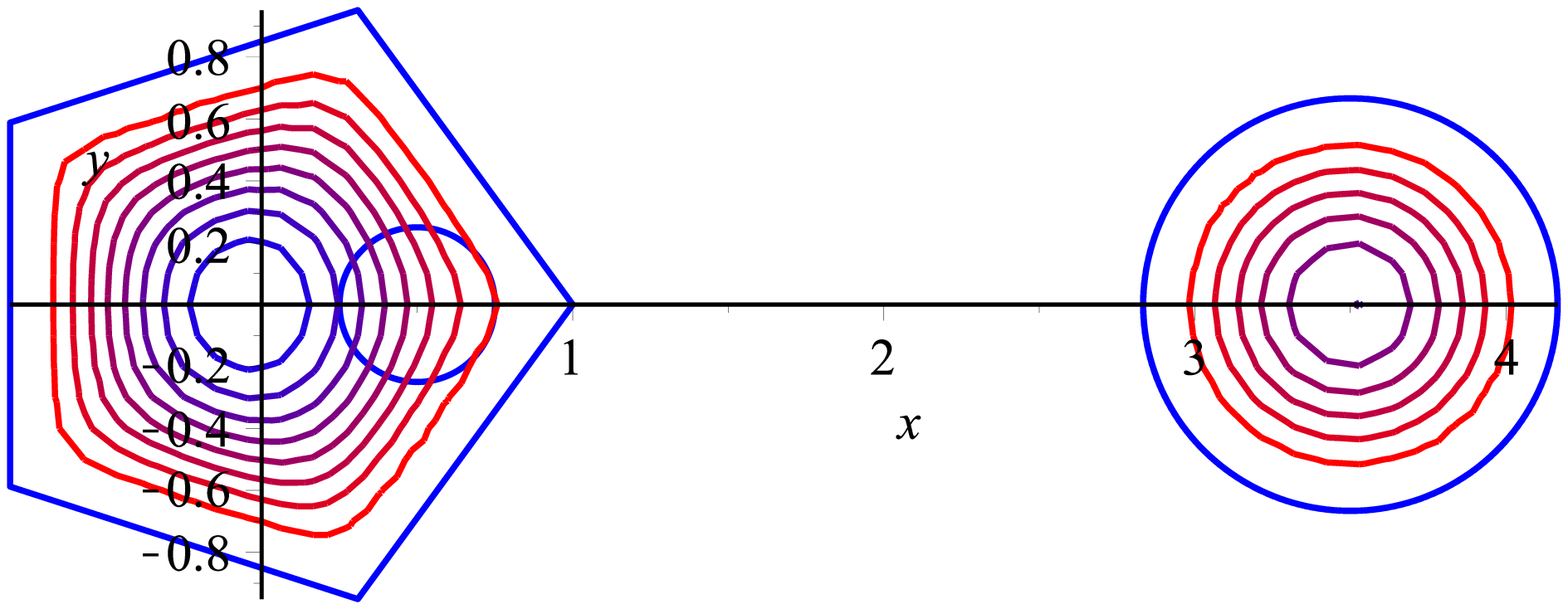}
\end{center}
\caption{\texttt{Phase A:} Level curves of $\lambda_{80}(G^\star,x+iy)$, on $\{(x,y):-2\le
x\le 5,-2\le y\le 2\}$, with $G$ as in Example~\ref{ex:diskpen}.}
\label{fig:diskpen-vili}
\end{figure}

In Figures~\ref{fig:diskpenA} and ~\ref{fig:diskpenB} we show the application of the two phases of the algorithm
on a frame that was suggested by the position of the zeros in Figure~\ref{fig:ro-diskpent80}.
In order to emphasize the importance of the information about zeros, we depict in Figure~\ref{fig:rec_disk-pent-disk80b-e}
the application of \texttt{Phase A}, with an arbitrarily chosen frame.

\begin{example}\label{ex:3disks-3disks}
Recovery for the archipelago of the three disks
$G_1=\{z:|z+1|<1/2\}$, $G_2=\{z:|z-2|<1\}$ and $G_3=\{z:|z-2i|<1/2\}$ and lake
$\mathcal{K}:=\cup_{j=1}^3\mathcal{K}_j$, where $\mathcal{K}_j$ are the following closed disks
$\mathcal{K}_1=\{z:|z+1|\le 1/3\}$, $\mathcal{K}_2=\{z:|z-2|\le 1/3\}$
and $\mathcal{K}_3=\{z:|z-2i|\le 1/4\}$.
\end{example}

\begin{figure}[h]
\begin{center}
\includegraphics*[scale=0.4]{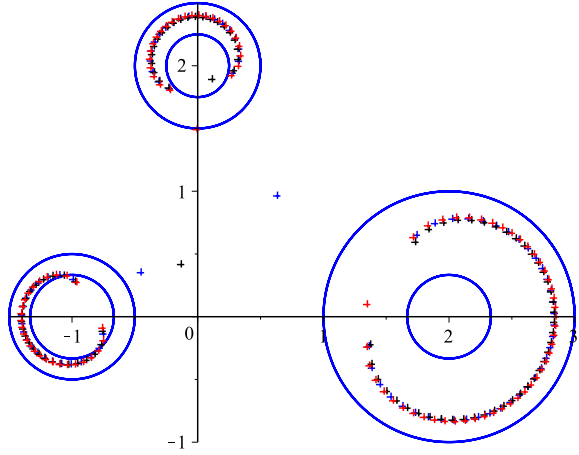}
\end{center}
\caption{Zeros of the polynomials $p_n(G^\star,z)$ of Example~\ref{ex:3disks-3disks}, for $n=80,90$ and $100$.}
\label{fig:ro-3disks}
\end{figure}

\begin{figure}[h]
\begin{center}
\includegraphics*[scale=0.4]{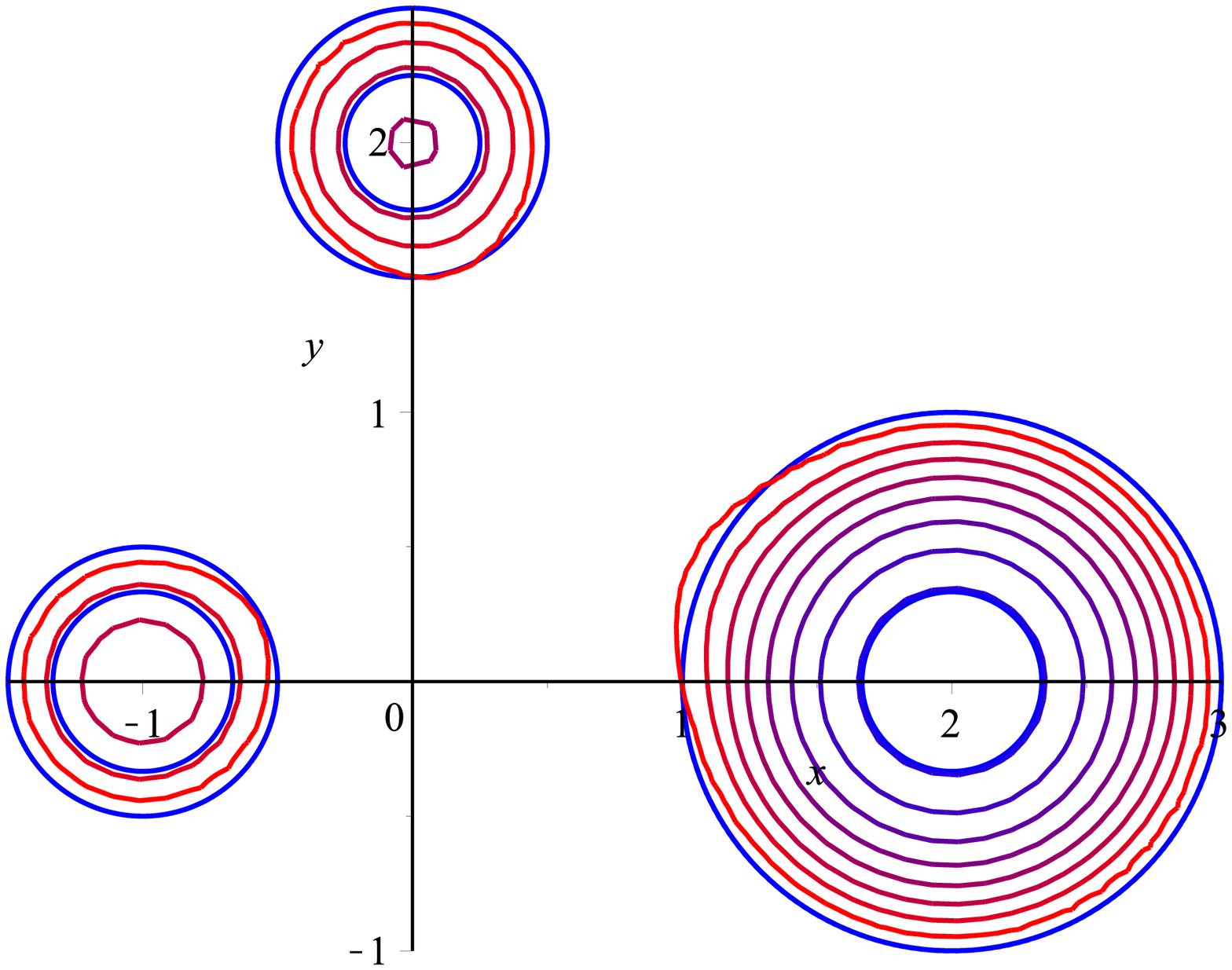}
\end{center}
\caption{\texttt{Phase A:} Level curves of ${\lambda^{1/2}_{100}}(G^\star,x+iy)$, on $\{(x,y):-3\le
x\le 4,-2\le y\le 3\}$, with $G^\star$ as in Example~\ref{ex:3disks-3disks}.}
\label{fig:rec-3disksA}
\end{figure}

\begin{figure}[h]
\begin{center}
\includegraphics*[scale=0.4]{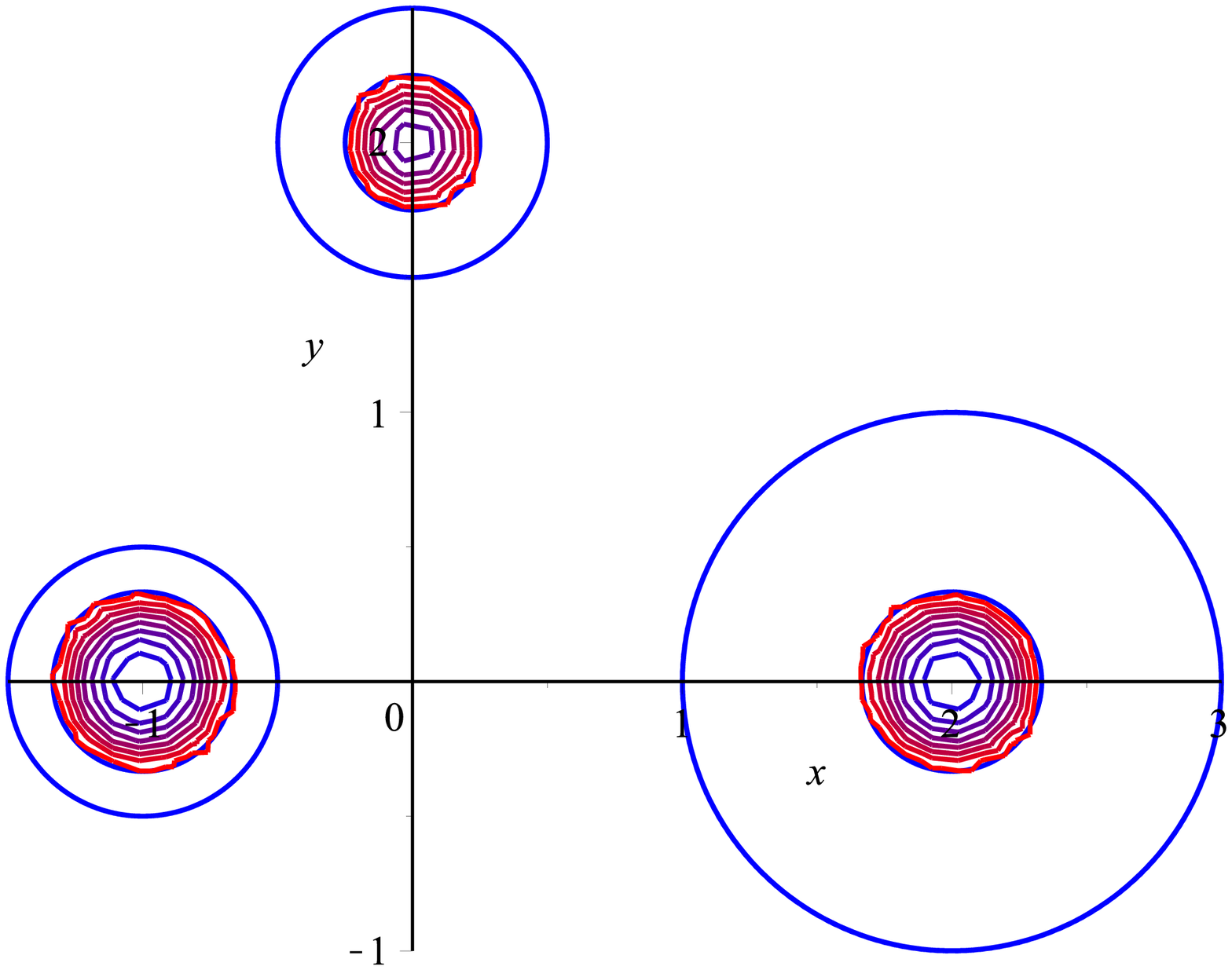}
\end{center}
\caption{\texttt{Phase B:} Level curves of ${\lambda^{1/2}_{100}}(\widehat{\mathcal{K}},x+iy)$, on $\{(x,y):-3\le
x\le 4,-2\le y\le 3\}$, with $G^\star$ as in Example~\ref{ex:3disks-3disks}.}
\label{fig:rec-3disksB}
\end{figure}

In Figure~\ref{fig:ro-3disks} we show the zeros $p_n(G^\star,z)$, for $n=80,90$ and $100$.
This should be compared with Figure 13 in \cite{GPSS}, which depicts zeros of $p_n(G,z)$.
In Figures~\ref{fig:rec-3disksA} and ~\ref{fig:rec-3disksB} we show the application of the two phases of the algorithm
on a frame that was suggested by the position of zeros in Figure~\ref{fig:ro-3disks}.

All the computations  were carried out on a MacBook Pro 2.4GHz Intel Core i7, using Maple 16.

\bigskip
\noindent E.B. Saff\\
Center for Constructive Approximation,\\
Department of Mathematics\\
Vanderbilt University\\
1326 Stevenson Center\\
37240 Nashville, TN\\
USA\\
{\it edward.b.saff@vanderbilt.edu}
\bigskip

\noindent Herbert Stahl\dag
\bigskip

\noindent N. Stylianopoulos\\
Department of Mathematics and Statistics,\\
University of Cyprus,\\
P.O. Box 20537, \\
1678 Nicosia, \\
Cyprus\\
{\it nikos@ucy.ac.cy}
\bigskip

\noindent Vilmos Totik\\
Bolyai Institute\\
MTA-SZTE Analysis and Stochastics Research Group\\
University of Szeged\\
Aradi v. tere 1, \\
6720  Szeged, \\
Hungary
\smallskip

\noindent and
\smallskip

\noindent Department of Mathematics and Statistics\\
University of South Florida \\
4202 E. Fowler Ave. CMC342\\
Tampa, FL, 33620  \\
USA\\
{\it totik@mail.usf.edu}


\clearpage
\begin{thebibliography}{10}
\bibitem{DuSc04}
P.~Duren and A.~Schuster, \emph{Bergman {S}paces}, Mathematical Surveys and
  Monographs, vol. 100, American Mathematical Society, Providence, RI, 2004.

\bibitem{Gabook87} D. Gaier,  {\it Lectures on {C}omplex {A}pproximation},
Birkh\"auser Boston Inc., Boston, MA, 1987.

\bibitem{GPSS}
B.~Gustafsson, M.~Putinar, E.~Saff, and N.~Stylianopoulos, \emph{Bergman
  polynomials on an archipelago: Estimates, zeros and shape reconstruction},
  Advances in Math. \textbf{222} (2009), 1405--1460.


\bibitem{PW86}
N.~Papamichael and M.~K. Warby, \emph{Stability and convergence properties of
  {B}ergman kernel methods for numerical conformal mapping}, Numer. Math.
  \textbf{48} (1986), no.~6, 639--669.

\bibitem{Pommerenke} Ch. Pommerenke, {\it Boundary Behavior of Conformal
Mappings}, Grundlehren der mathematischen Wissenschaften,
{\bf 299}, Springer Verlag, Berlin, Heidelberg New York, 1992.



\bibitem{Saff1} E.~B. Saff,
\emph{Remarks on relative asymptotics for general orthogonal
  polynomials}, Recent trends in orthogonal polynomials and approximation
  theory, Contemp. Math., vol. 507, Amer. Math. Soc., Providence, RI, 2010,
  pp.~233--239.

\bibitem{SafSty} E.B. Saff and N. Stylianopoulos, \emph{ Asymptotics for Hessenberg Matrices for the Bergman Shift Operator on Jordan Regions}, Complex Anal. Oper. Theory \textbf{8} (2014), no. 1, 1–-24.

\bibitem{SaffTotik}
E.~B. Saff and V.~Totik, \emph{Logarithmic {P}otentials with {E}xternal {F}ields},
 Grundlehren der mathematischen Wissenschaften,
{\bf 316}, Springer-Verlag, New York/Berlin, 1997.

\bibitem{Si12}
B.~Simanek, \emph{A new approach to ratio asymptotics for orthogonal
  polynomials}, J. Spectr. Theory \textbf{2} (2012), no.~4, 373--395.

\bibitem{StahlTotik} H. Stahl and V. Totik, {\it General Orthogonal Polynomials},
Encyclopedia of Mathematics, {\bf 43},
Cambridge University Press, New York 1992

\bibitem{St-CA13}
N.~Stylianopoulos, \emph{Strong asymptotics for {B}ergman polynomials over
domains with corners and applications}, Constr. Approx. \textbf{38} (2013), 59--100.

\bibitem{Su74} P. K. Suetin, {\it Polynomials orthogonal over a region and Bieberbach
              polynomials}, Proceedings of the Steklov Institute of Mathematics, {\bf 100} (1971),
              Amer. Math. Soc., Providence, R.I., 1974

\bibitem{To05}
V.~Totik, \emph{Orthogonal polynomials}, Surv. Approx. Theory \textbf{1}
  (2005), 70--125 (electronic).

\bibitem{Walsh} J. L. Walsh,
{\it Interpolation and Approximation by
Rational Functions in the Complex Domain}, fifth edition,
Amer. Math. Soc. Colloquium Publications, {\bf XX},
Amer. Math. Soc., Providence, 1969.
\end{thebibliography}
\end{document}